\documentclass[a4paper,12pt]{article}
\usepackage{amssymb}
\usepackage{amsmath}
\usepackage{amsthm}
\usepackage{graphicx}
\usepackage{hyperref}
\usepackage[square,numbers,sort&compress]{natbib}
\newtheorem{thm}{Theorem}[section]
\newtheorem{cor}{Corollary}[section]
\newtheorem{lem}{Lemma}[section]

\theoremstyle{definition}

\theoremstyle{remark}
\newtheorem{rem}{Remark}[section]

\numberwithin{equation}{section}

\begin{document}
\begin{center}
{\Large \bf Some compact generalization of Bernstein-type inequalities preserved by modified Smirnov operator}
\end{center}
\begin{center}
		{\normalsize \bf Deepak Kumar\footnote{deepak.kumar26@s.amity.edu}, D. Tripathi \footnote{dinesh@mru.edu.in, dn3pathi.phdcpa@gmail.com} and Sunil Hans\footnote{ corresponding author (sunil.hans82@yahoo.com)}}
\end{center}
\begin{center}
 {\normalsize $^{1 , 3}$Department of Applied Mathematics, Amity University,\\ Noida-201313, India.}\\
{\normalsize $^{2}$Department of Science - Mathematics, School of Sciences\\ Manav Rachana University, Faridabad-121004, India.}\\
\end{center}
 \begin{abstract}  
 Let $P(z)$ be a polynomial of degree $n$. In $2004$, Aziz and Rather \cite{aziz2004some} investigated the dependence of  \[\bigg|P(Rz)-\alpha P(z)+\beta\biggl\{\biggl(\frac{R+1}{2}\biggr)^n-|\alpha|\biggr\}P(z)\bigg|, \ \text{for} \ z \in B(\mathbb{D}),\] on $\max_{z\in B(\mathbb{D})}|P(z)|$, for every real and complex number $\alpha, \beta$ satisfying $|\alpha| \leq 1$, $|\beta| \leq 1$, and $R \geq 1$. This paper presents a compact generalization of several well-known polynomial inequalities using modified Smirnov operator, demonstrating that the operator preserves inequalities between polynomials.  
 \end{abstract}
 \begin{center}
     \section{INTRODUCTION}
     \end{center}
     Let $P(z)$ be a polynomial of degree $n$ and $P'(z)$ be the derivative of polynomial $P(z)$. Let $\mathbb{D}$ be the open unit disk $\{z\in\mathbb{C}; |z|<1\} $, so that $\overline{\mathbb{D}}$ is it's closure and $B(\mathbb{D) }$ denotes its boundary, then
     \begin{align}{\label{1.1}}
         \max_{z\in B(\mathbb{D})}|P'(z)| \leq n \max_{z\in B(\mathbb{D})}|P(z)|,
     \end{align}
     and
     \begin{align}{\label{1.2}}
         \max_{z\in B(\mathbb{D})}|P(Rz)| \leq R^n \max_{z\in B(\mathbb{D})}|P(z)|.
     \end{align}
    Inequality (\ref{1.1}) can be obtained by a direct result of a theorem of S. Bernstein \cite{bernstein1912ordre} for the derivative of a polynomial. A straightforward inference from the maximum modulus principle \cite{polya1925aufgaben} yields the inequality (\ref{1.2}). In both inequalities, the equality holds for $P(z)=\lambda z^n, \lambda\neq 0$. Inequalities (\ref{1.1}) and (\ref{1.2}) can be improved if the zeros are restricted. Erd\"os conjectured and Lax \cite{lax1944proof} proved that if $P(z)$ has no zeros in $\mathbb{D}$, then
\begin{align}{\label{1.3}}
        \max_{z\in B(\mathbb{D})}|P'(z)| \leq \frac{n}{2}\max_{z\in B(\mathbb{D})}|P(z)|,
    \end{align}
    and for $R \geq 1$,
    \begin{align}{\label{1.4}}
        \max_{z\in B(\mathbb{D})}|P(Rz)| \leq \frac{R^n+1}{2}\max_{z\in B(\mathbb{D})}|P(z)|.
    \end{align}
   Ankeny and Rivlin \cite{ankeny1955theorem}  used (\ref{1.3}) to prove the above inequality (\ref{1.4}). As a generalization of inequalities  (\ref{1.1}) and (\ref{1.2}), Aziz and Rather \cite{aziz1999inequality} proved that if $P(z)$ is a polynomial of degree $n$, then for every real and complex number $\alpha$ with $|\alpha|\leq 1$ and $R\geq 1$,
   \begin{align}{\label{1.5}}
       |P(Rz)-\alpha P(z)|\leq |R^n-\alpha||z|^n\max_{z \in B(\mathbb{D})}|P(z)| \ \text{for} \ z \in \mathbb{C\backslash D}.
   \end{align}
   The above result is best possible and equality in (\ref{1.5}) holds for $P(z)=\lambda z^n, \lambda\neq 0$. Inequality (\ref{1.1}) follows from inequality (\ref{1.5}) on dividing it by $R-1$ and taking the limit as $ R \to 1 $ with $ \alpha = 1 $. Inequality (\ref{1.5}) reduces to (\ref{1.2}) for $\alpha = 0$.
      As an improvement of inequality (\ref{1.5}), the authors \cite{aziz1999inequality} have also shown that if $P(z)\neq 0$ in $\mathbb{D}$, then for every real and complex number $\alpha$ with $|\alpha|\leq 1$ and $R\geq 1$
   \begin{align}{\label{1.6}}
        |P(Rz)-\alpha P(z)|\leq \frac{1}{2}\{|R^n-\alpha||z|^n+|1-\alpha|\}\max_{z \in B(\mathbb{D})}|P(z)|,
   \end{align}
   for $ z \in \mathbb{C\backslash D}$.
   The result is sharp and equality in (\ref{1.6}) holds for $ P(z)=z^n + 1$. After dividing both sides of (\ref{1.6}) by $(R-1)$ and taking limit $R \to 1$ with $\alpha=1$, inequality (\ref{1.3}) has been obtained. Furthermore, inequality (\ref{1.6}) reduces to (\ref{1.4}) for $\alpha =0$. \\
   In $1930$, S. Bernstein \cite{bernstein1930limitation} also proved the following result:
   \begin{thm}{\label{thmA}}
       Let $F(z)$ be a polynomial of degree $n$, having all its zeros in $\overline{\mathbb{D}}$ and $P(z)$ be a polynomial of degree not exceeding that of $F(z)$. If $|P(z)|\leq |F(z)|$  on $B(\mathbb{D})$, then 
       \begin{align*}
          |P'(z)| \leq |F'(z)| \ for \ z\in \mathbb{C} \backslash \mathbb{D}.
       \end{align*} 
       The equality holds only if $P(z)=e^{i\gamma}F(z), \gamma \in \mathbb{R}$.
   \end{thm}
   For $z\in \mathbb{C}\backslash \mathbb{D}$, denoting $\Omega _{|z|}$ the image of the disk $\{t\in \mathbb{C}; |t|<|z| \}$ under the mapping $\phi(t)= \frac{t}{t+1}$, Smirnov \cite{smirnov1964constructive} as a generalization of Theorem \ref{thmA} proved the following:
   \begin{thm}{\label{B}}
      Let $F$ and $P$ be the polynomial of possessing condition as in Theorem \ref{thmA}. Then  for $z\in \mathbb{C\backslash D}$
      \begin{align}{\label{1.7}}
          |\mathbb{S}_\alpha[P](z)| \leq |\mathbb{S}_\alpha[F](z)|,
            \end{align}
          for all $\alpha \in \overline{\Omega}_{|z|}$, with $\mathbb{S}_\alpha[P](z)= z P'(z)-n\alpha P(z)$, where $\alpha$ is a constant.  
   \end{thm}
   For $\alpha \in \overline{\Omega}_{|z|}$ in the inequality (\ref{1.7}), the equality is holds at a point $z \in \mathbb{C\backslash D}$, only if $P(z) = e^{i \gamma}F(z), \gamma \in \mathbb{R}$. We note that for fixed $z \in \mathbb{C \backslash D}$, the inequality (\ref{1.7}) can be replaced by (see for reference \cite{ganenkova2019variations})
   \[\biggl|zP'(z)-n\frac{az}{1+az}P(z)\biggr| \leq \biggl|zF'(z)-n\frac{az}{1+az}F(z)\biggr|,\]
   where $a$ is an arbitrary number in $\overline{\mathbb{D}}$.\\
   Equivalently for $z \in \mathbb{C\backslash D}$
   \[|\Tilde{\mathbb{S}}_a[P](z)| \leq |\Tilde{\mathbb{S}}_a[F](z)|,\]
   where $\Tilde{\mathbb{S}}_a[P](z)=(1+az)P'(z)-naP(z)$ is known as the modified Smirnov operator. The modified Smirnov operator $\Tilde{\mathbb{S}}_a$ is preferred over the Smirnov operator $\mathbb{S}_\alpha$ because its parameter  $a$ is independent of $z$, while the parameter $\alpha$ in $\mathbb{S}_\alpha$ depends on $ z$.\\
   Shah and Fatima \cite{shah2022bernstein} used modified Smirnov operator to generalize inequalities (\ref{1.1}), (\ref{1.2}), (\ref{1.3}) and (\ref{1.4}) and proved that if $P(z)$ is a polynomial of degree $n$, such that $|P(z)| \leq M$ for $z\in B(\mathbb{D})$, then 
 for $ z \in \mathbb{C\backslash D}$
   \begin{align}{\label{1.8}}
       |\tilde{\mathbb{S}}_a[P](z)|\leq M|\tilde{\mathbb{S}}_a[z^n]|,
       \end{align}
       equivalently
       \begin{align}{\label{1.9}}
           |(1+az)P'(z)-naP(z)| \leq Mn|z|^{n-1},
       \end{align}
and if $P(z)\neq 0 $ in $\mathbb{D}$, then for $ z \in\mathbb{C\backslash D}$
   \begin{align}{\label{1.10}}
    |\tilde{\mathbb{S}}_a[P](z)|\leq \frac{1}{2}\{|\tilde{\mathbb{S}}_a[z^n]|+n|a|\}\max_{z\in B(\mathbb{D})}|P(z)|,
   \end{align}
   equivalently
   \begin{align}\label{1.11}
       |(1+az)P'(z)-naP(z)| \leq \frac{1}{2}\{n|z|^{n-1}+n|a|\}\max_{z \in B(\mathbb{D})}|P(z)|.
   \end{align}
   By setting $a=0$, inequalities $(\ref{1.9})$ and (\ref{1.11}) reduce to inequalities (\ref{1.1}) and (\ref{1.3}) respectively. Likewise, by setting $ a = -\frac{1}{z}$ with $z = Re^{i\theta}, R \geq 1$, inequalities (\ref{1.2}) and (\ref{1.4}) can be obtained from inequalities (\ref{1.9}) and (\ref{1.11}), respectively.\\
   Wani and Liman \cite{wani2024bernstein} have generalized inequality (\ref{1.5}) concerning the modified Smirnov operator and proved that if $P(z)$ is a polynomial of degree $n$, then for every real and complex number $\alpha$ with $|\alpha| \leq 1$ and $R \geq 1$ 
\begin{align}{\label{1.12}}
       |\tilde{\mathbb{S}}_a[P](Rz)-\alpha \tilde{\mathbb{S}}_a[P](z)| \leq |R^n-\alpha||\tilde{\mathbb{S}}_a[z^n]|\max_{z\in B(\mathbb{D})}|P(z)|,
   \end{align}
    for $z \in \mathbb{C\backslash D}$. The result is sharp and equality holds in (\ref{1.12}) for $ P(z)=\lambda z^n, \lambda \neq 0$.
   Further, as a generalization of inequality (\ref{1.6}), the authors \cite{wani2024bernstein} have also shown that if $P(z) \neq 0$ in $\mathbb{D}$, then for every real and complex number $\alpha$ with $|\alpha|\leq 1$ and $R \geq 1$
   \begin{align}{\label{1.13}}
       |\tilde{\mathbb{S}}_a[P](Rz)-\alpha\tilde{\mathbb{S}}_a[P](z)| \leq \biggl\{\frac{|R^n-\alpha||\tilde{\mathbb{S}}_a[z^n]|+n|1-\alpha||a|}{2}\biggr\}\max_{z \in B(\mathbb{D})}|P(z)|, 
   \end{align}
    for $\ z\in \mathbb{C\backslash D}$. The result is the best possible and the equality holds for $P(z)= \lambda z^n, \lambda \neq 0$. 
   Aziz and Rather \cite{aziz2004some} have investigated the dependence of
   \[\bigg|P(Rz)-\alpha P(z)+\beta\biggl\{\biggl(\frac{R+1}{2}\biggr)^n-|\alpha|\biggr\}P(z)\bigg| \ \text{for} \ z \in B(\mathbb{D}) \ \text{on} \ \max_{z \in B(\mathbb{D})}|P(z)|\]
   for every real and complex number $\alpha,\beta$ with $|\alpha|\leq 1, |\beta|\leq1$ and $R\geq 1$. As a compact generalization of the inequalities (\ref{1.1}), (\ref{1.2}) and (\ref{1.5}), they provide the following theorem:
   \begin{thm}{\label{thmC}}
       If P(z) is a polynomial of degree n, then for every real and complex number $ \alpha,\beta$ with $|\alpha| \leq 1, |\beta| \leq 1 \  and \ R \geq 1$,
       \begin{align}{\label{1.14}}
          \nonumber &\bigg|P(Rz)-\alpha P(z)+\beta\biggl\{\biggl(\frac{R+1}{2}\biggr)^n-|\alpha|\biggr\}P(z)\bigg|\\& \leq \bigg|R^n-\alpha+\beta\biggl\{\biggl(\frac{R+1}{2}\biggr)^n-|\alpha|\biggr\}\bigg||z|^n\max_{z \in B(\mathbb{D})}|P(z)| \ for \ z \in \mathbb{C\backslash D}.
       \end{align}
        The results is sharp and the equality in \textnormal{(\ref{1.14})} holds for $P(z)=\lambda z^n, \lambda \neq 0$.
   \end{thm}
    As an improvement of the above result, Aziz and Rather \cite{aziz2004some} proved the following theorem for the class of polynomials having no zero in the unit disk.
    \begin{thm}{\label{thmD}}
        If $P(z)$ is a polynomial of degree $n$, which does not vanish in $\mathbb{D}$, then for every real and complex number $\alpha, \beta$ with $|\alpha| \leq 1, |\beta| \leq 1$ and $R\geq 1$,
        \begin{align}{\label{1.15}}
            \nonumber &\bigg|P(Rz)-\alpha P(z)+\beta\biggl\{\biggl(\frac{R+1}{2}\biggr)^n-|\alpha|\biggr\}P(z)\bigg|\\& \nonumber\leq \frac{1}{2}\biggl[\bigg|R^n-\alpha+\beta\biggl\{\biggl(\frac{R+1}{2}\biggr)^n-|\alpha|\biggr\}\bigg||z|^n\\& +\bigg|1-\alpha + \beta\biggl\{\bigg(\frac{R+1}{2}\bigg)^n-|\alpha|\biggr\}\bigg|\bigg]\max_{z \in B(\mathbb{D})}|P(z)|.
        \end{align}
        The result is best possible and equality in \textnormal{(\ref{1.15})} holds for $P(z)= z^n+1$.
    \end{thm}
    The above theorem presents a generalization of inequalities (\ref{1.3}), (\ref{1.4}), and (\ref{1.6}).\\
    In this paper, we prove the following results, which generalize Theorem \ref{thmC} and Theorem \ref{thmD} for the modified Smirnov operator.
 \begin{center}
     \section{MAIN RESULT}
 \end{center}
     \begin{thm}{\label{thm1}}
         If $P(z)$ is a polynomial of degree $n$, then for every real and complex number $\alpha, \beta$ with $|\alpha|\leq 1, |\beta|\leq 1$ and $R\geq 1$,
         \begin{align}{\label{2.1}}
            \nonumber &\biggl| \Tilde{\mathbb{S}}_a[P](Rz)-\alpha \Tilde{\mathbb{S}}_a[P](z)+\beta \biggl\{\biggl(\frac{R+1}{2}\biggr)^n-|\alpha| \biggr\}\Tilde{\mathbb{S}}_a[P](z)\biggr|\\& \leq \biggl|R^n-\alpha +\beta\biggl\{\biggl(\frac{R+1}{2}\biggr)^n -|\alpha|\biggr\}\biggr| |\Tilde{\mathbb{S}}_a[z^n]|\max_{z\in B(\mathbb{D})}|P(z)|,
         \end{align}
        for $\  z \in \mathbb{C \backslash D}$. The result is sharp and equality in \textnormal{(\ref{2.1})} holds for $P(z)=\lambda z^n, \lambda\neq 0.$
     \end{thm}
     \begin{rem}
        If in inequality (\ref{2.1}), we take $a=0$, we get the following result proved by Baseri et al. \cite{41bf3c67-6158-38a2-a5b2-f9e1437abf30}
        \begin{align*}
        \biggl|& RP'(Rz)-\alpha P'(z)+ \beta\biggl\{\biggl(\frac{R+1}{2}\biggr)^n -|\alpha|\biggr\}P'(z)\bigg|\\&  \leq \bigg|R^n-\alpha+\beta \biggl\{ \biggl( \frac{R+1}{2}\biggr)^n-|\alpha|\biggr\}\bigg|n|z|^{n-1}\max_{z \in B(\mathbb{D})}|P(z)|, \ for \ z\in \mathbb{C\backslash D}
         \end{align*}
        Inequality (\ref{1.14}) is a special case of inequality (\ref{2.1}) for $ a=-\frac{1}{z}$.
     \end{rem}
     \begin{rem}
      If we choose $\beta =0$ in (\ref{2.1}), then we get inequality (\ref{1.12}). In addition, if we consider $\alpha$ and $\beta$ to be zero, then we get
      \begin{align*}
      |\tilde{\mathbb{S}}_a[P](Rz)|\leq R^n|\tilde{\mathbb{S}}_a[z^n]|\max_{z \in B(\mathbb{D})}|P(z)|, \ for \ z\in \mathbb{C\backslash D}.
      \end{align*}
      By substituting, $a = -\frac{1}{z}$ and $R= 1$ in the above inequality, we obtain inequalities (\ref{1.2}) and (\ref{1.8}), respectively.
     \end{rem}
     If we choose $\alpha =1$ in inequality (\ref{2.1}) and divide it by $R-1$ and taking $ R \to 1$, we get the following result.
     \begin{cor}{\label{cor2.1}}
         If $P(z)$ is a polynomial of degree $n$, then for every real and complex number $ \beta$ with $|\beta|\leq 1$ and $R\geq 1$,
         \begin{align}{\label{2.2}}
            \bigg|z\tilde{\mathbb{S}}_a[P'](z)+\frac{n}{2}  \beta  
  \tilde{\mathbb{S}}_a[P](z) +P'(z)\biggr|\leq n\bigg|1+\frac{\beta}{2}\bigg||\tilde{\mathbb{S}}_a[z^n]|\max_{z \in B(\mathbb{D})}|P(z)|. 
         \end{align}
         for $z\in \mathbb{C\backslash D}$. The equality holds in \textnormal{(\ref{2.2})} for $  P(z)=\lambda z^n, \lambda \neq 0$.       
     \end{cor}
      For the $ \alpha =0$ inequality (\ref{2.1}), yield the following result, which generalizes one of the results of Jain \cite{jain1992maximum} for the modified Smirnov operator.
     \begin{cor}
     If $P(z)$ is a polynomial of degree $n$, then for every real complex number $ \beta$ with $|\beta|\leq 1$ and $R\geq 1$,
     \begin{align*}
         \nonumber \biggl| \Tilde{\mathbb{S}}_a[P](Rz)+\beta \biggl(\frac{R+1}{2}\biggr)^n\Tilde{\mathbb{S}}_a[P](z)\biggr| \leq  \biggl|R^n +\beta\biggl(&\frac{R+1}{2}\biggr)^n \biggr| |\Tilde{\mathbb{S}}_a[z^n]|\max_{z\in B(\mathbb{D})}|P(z)|,
         \end{align*}
          for  $z\in \mathbb{C \backslash D}$.\\ The above result is sharp and equality holds for $ P(z)=\lambda z^n, \lambda \neq 0.$
     \end{cor} 
     \begin{thm}{\label{thm2}}
     If $P(z)$ is a polynomial of degree $n$ and $Q(z)=z^n \overline{P(\frac{1}{\overline{z}})}$, then for every real and complex number $\alpha, \beta$ with $|\alpha|\leq1,|\beta|\leq 1$ and $R\geq 1$,
     \begin{align}{\label{2.3}}
         \nonumber& \biggl|\Tilde{\mathbb{S}}_a[P](Rz)-\alpha\Tilde{\mathbb{S}}_a[P](z)+\beta\biggl\{\biggl(\frac{R+1}{2}\biggr)^n-|\alpha|\biggr\}\Tilde{\mathbb{S}}_a[P](z)\biggr|+\\& \nonumber \biggl|\Tilde{\mathbb{S}}_a[Q](Rz)-\alpha\Tilde{\mathbb{S}}_a[Q](z)+\beta\biggl\{\bigg(\frac{R+1}{2}\biggr)^n-|\alpha|\biggr\}\Tilde{\mathbb{S}}_a[Q](z)\biggr|  \\& \nonumber \leq \biggl[ \biggl|R^n-\alpha+\beta\biggl\{\bigg(\frac{R+1}{2}\bigg)^n-|\alpha|\biggr\}\biggr||\Tilde{\mathbb{S}}_a[z^n]|+ \\&\qquad \biggl|1-\alpha +\beta\bigg\{\biggl(\frac{R+1}{2}\biggr)^n-|\alpha|\bigg\}\biggr|n|a|\biggr]\max_{z\in B(\mathbb{D})}|P(z)|,\ for \ z \in \mathbb{C\backslash D}
        \end{align}
        The result is sharp and equality holds for $P(z)= \lambda z^n,  \lambda \neq 0$.
     \end{thm}
     If we take $\alpha=0$ in the inequality (\ref{2.3}), the following result, which extends a result is given by Jain \cite{jain1992maximum} for the modified Smirnov operator.
     \begin{cor}
          If $P(z)$ is a polynomial of degree $n$, then for every real and complex number $\beta$ with $|\beta|\leq 1$ and $R\geq 1$
          \begin{align*}
              \nonumber& \biggl|\Tilde{\mathbb{S}}_a[P](Rz)+\beta\biggl(\frac{R+1}{2}\biggr)^n\Tilde{\mathbb{S}}_a[P](z)\biggr|+ \biggl|\Tilde{\mathbb{S}}_a[Q](Rz)+\beta\bigg(\frac{R+1}{2}\biggr)^n\Tilde{\mathbb{S}}_a[Q](z)\biggr|  \\& \nonumber \leq \biggl[ \biggl|R^n +\beta\bigg(\frac{R+1}{2}\bigg)^n\biggr||\Tilde{\mathbb{S}}_a[z^n]|+ \biggl|1 +\beta \biggl(\frac{R+1}{2}\biggr)^n\biggr|n|a|\biggr]\max_{z\in B(\mathbb{D})}|P(z)|,
          \end{align*}
          for $z \in \mathbb{C\backslash D}$ and $Q(z)=z^n \overline{P\left(\frac{1}{\overline{z}}\right)}$. The result is sharp and equality holds for $P(z)=\lambda z^n, \lambda \neq 0$.
     \end{cor}
      In inequality (\ref{2.3}), if we consider $\alpha=1$ and divide both sides by $R-1$ and letting limit $R \to 1$, we have  
     \begin{cor}
         If $P(z)$ is a polynomial of degree $n$ and $Q(z)=z^n \overline{P(\frac{1}{\overline{z}})}$, then for every real and complex number $ \beta$ with $|\beta|\leq 1$, $R\geq 1$ and $z \in \mathbb{C\backslash D}$
         \begin{align}{\label{2.4}}
             \nonumber& \biggl|\Tilde{\mathbb{S}}_a[P](Rz)+\beta\biggl(\frac{R+1}{2}\biggr)^n\Tilde{\mathbb{S}}_a[P](z) +P'(z)\biggr|+\\& \nonumber \biggl|\Tilde{\mathbb{S}}_a[Q](Rz)+\beta\bigg(\frac{R+1}{2}\biggr)^n\Tilde{\mathbb{S}}_a[Q](z) +Q'(z)\biggr|  \\& \leq n\biggl[\bigg|1+\frac{\beta}{2}\bigg||\tilde{\mathbb{S}}_a[z^n]|+\frac{n}{2}|\beta||a|\biggl]\max_{z \in B(D)}|P(z)|, \ for \ z \in \mathbb{C \backslash D}
         \end{align}
         The above result is best possible and equality holds for $P(z)= \lambda z^n, \lambda\neq 0$.
     \end{cor}
     \begin{rem}
         For $\beta =0$, Theorem \ref{thm2}, yield a result of Wani and Liman \cite{wani2024bernstein} and a result of Shah and Fatima \cite{shah2022bernstein} follows from it, when $\alpha =\beta=0$ and $R=1$.
     \end{rem}
     \begin{rem}
         If we choose $a=0$ in inequality (\ref{2.3}), we get
         \begin{align*}
         &\biggl[RP'(Rz)-\alpha P'(z)+\beta\biggl\{\bigg(\frac{R+1}{2}\bigg)^n-|\alpha|\biggr\}P'(z)\biggr]+\\ &\biggl[RQ'(Rz)-\alpha Q'(z)+\beta\biggl\{\bigg(\frac{R+1}{2}\bigg)^n-|\alpha|\biggr\}Q'(z)\biggr] \\ &\leq \bigg|R^n-\alpha+\beta \biggl\{ \biggl( \frac{R+1}{2}\biggr)^n-|\alpha|\biggr\}\bigg|n|z|^{n-1}\max_{z \in B(\mathbb{D})}|P(z)|, \ for \ z\in \mathbb{C\backslash D} 
         \end{align*}
         and for $a=-\frac{1}{z}$, Theorem \ref{thm2} reduce to a result of Aziz and Rather \cite{aziz2004some}.
     \end{rem}
     \begin{thm}{\label{thm3}}
         If $P(z)$ is a polynomial of degree $n$, which doesn't vanish in $\mathbb{D}$, then for every real and  complex number $\alpha, \beta$ with $|\alpha|\leq1,|\beta|\leq 1$ and $R\geq 1$,
         \begin{align}{\label{2.5}}
           \nonumber  \biggl|&\Tilde{\mathbb{S}}_a[P](Rz)-\alpha \Tilde{\mathbb{S}}_a[P](z)+\beta \biggl\{\biggl(\frac{R+1}{2}\biggr)^n -|\alpha|\biggr\}\Tilde{\mathbb{S}}_a[P](z)\biggr|\\& \nonumber\leq \frac{1}{2}\biggl[ \biggl|R^n-\alpha+\beta\biggl\{\bigg(\frac{R+1}{2}\bigg)^n-|\alpha|\biggr\}\biggr||\Tilde{\mathbb{S}}_a[z^n]|+\\& \qquad\biggl|1-\alpha + \beta\bigg\{\biggl(\frac{R+1}{2}\biggr)^n-|\alpha|\bigg\}\biggr|n|a|\biggr]\max_{z \in B(\mathbb{D})}|P(z)|\ for\ z\in \mathbb{C\backslash D}.
         \end{align}
         The equality in \textnormal{(\ref{2.5})} holds for $P(z)=z^n+1$.
     \end{thm}
     If we consider $\alpha = 0$ in inequality (\ref{2.5}), we derive the following result, which generalizes a particular result proved by Jain \cite{jain1997generalization}.
     \begin{cor}
          If $P(z)$ is a polynomial of degree $n$, which doesn't vanish in $\mathbb{D}$, then for every real and  complex number $ \beta $ with $|\beta|\leq 1$ and $R\geq 1$
     \begin{align*}
         \nonumber  \biggl|\Tilde{\mathbb{S}}_a[P](Rz)+\beta \biggl(\frac{R+1}{2}\biggr)^n &\Tilde{\mathbb{S}}_a[P](z)\biggr| \leq \frac{1}{2}\biggl[ \biggl|R^n +\beta\bigg(\frac{R+1}{2}\bigg)^n\biggr||\Tilde{\mathbb{S}}_a[z^n]|+\\& \biggl|1 +\beta\biggl(\frac{R+1}{2}\biggr)^n \biggr|n|a|\biggr]\max_{z \in B(\mathbb{D})}|P(z)| \ for\ z\in \mathbb{C\backslash D}.
     \end{align*}
     The above result is sharp and holds for $P(z)=z^n+1$.
     \end{cor}
    The next corollary is obtained by taking $\alpha=1$ in Theorem \ref{thm3}, dividing by $R-1$, and then letting $R\to 1$. This gives a refinement of Corollary \ref{cor2.1} for polynomials not vanishing in the unit disk.
     \begin{cor}
            If $P(z)$ is a polynomial of degree $n$, which doesn't vanish in $\mathbb{D}$, then for every real and complex number $ \beta $ with $|\beta|\leq 1$ and $R\geq 1$
            \begin{align}{\label{2.6}}
                \nonumber \bigg|&z\tilde{\mathbb{S}}_a[P'](z)+\frac{n}{2}  \beta   \tilde{\mathbb{S}}_a[P](z) +P'(z)\biggr| \\ &\leq \frac{n}{2}\biggl[\bigg|1+\frac{\beta}{2}\bigg||\tilde{\mathbb{S}}_a[z^n]|+\frac{n}{2}|\beta||a|\biggr]\max_{z \in B(\mathbb{D})}|P(z)|, \ for\ z\in \mathbb{C\backslash D}.
            \end{align}
            The result is the best possible and equality in \textnormal{(\ref{2.6})} holds for $P(z)=z^n+1$.
     \end{cor}
     \begin{rem}
   Theorem \ref{thm3} reduces to inequality (\ref{1.13}) for $\beta=0$. If we consider $\alpha=\beta=0 $ and $ R=1$ in (\ref{2.5}), it reduces to the inequality (\ref{1.10}). Further, if we choose $a=-\frac{1}{z}$ in (\ref{2.5}), we get inequality (\ref{1.15}).
     \end{rem}
 \begin{center}
\section{LEMMA}	
\end{center}
     For the proof of these theorems, we need the following lemmas.
     The first lemma is due to Aziz \cite{aziz1987growth}.
\begin{lem}{\label{lem1}}
 If $P(z)$ is a polynomial of degree $n$, having all its zeros in $|z| \leq k, k\leq 1$ , then every $R\geq 1$ 
 \begin{align}{\label{3.1}}
     |P(Rz)| \geq \left(\frac{R+k}{1+k}\right)^n|P(z)|, \ for\ z \in B(\mathbb{D}). 
 \end{align}
\end{lem}
\begin{lem}{\label{lem2}}
 Let $P(z)$ be a polynomial of degree $n$ with all its zeros in $\overline{\mathbb{D}}$. If $a\in B(\mathbb{D})$ is not the exceptional value for $P(z)$, then all the zeros of $\Tilde{\mathbb{S}}_a[P](z)$ lie in ${\overline{\mathbb{D}}}$.  
\end{lem}
The above lemma is due to Genenkova and Starkov \cite{ganenkova2019variations}. In addition, the next two lemmas are given by Shah and Fatima \cite{shah2022bernstein}. 
\begin{lem}{\label{lem3}}
   If $ P(z)$ is a polynomial of degree $n$, having no zeros in $\mathbb{D}$, then
   \begin{align}{\label{3.2}}
       |\Tilde{\mathbb{S}}_a[P](z)|\leq  |\Tilde{\mathbb{S}}_a[Q](z)|, \quad for \ z \in \mathbb{C\backslash D},
   \end{align}
   where $Q(z)=z^n \overline{P(\frac{1}{\overline{z}})}$.
   \end{lem}
   \begin{lem}{\label{lem4}}
   If $P(z)$ is a polynomial of degree $n$, then for \ $z \in \mathbb{ C\backslash D}$
   \begin{align}{\label{3.3}}
    |\Tilde{\mathbb{S}}_a[P](z)| +  |\Tilde{\mathbb{S}}_a[Q](z)| \leq \{|\Tilde{\mathbb{S}}_a[z^n]| + n|a|\} \max_{z \in B(\mathbb{D})}|P(z)|,
    \end{align}
     where $Q(z)=z^n \overline{P(\frac{1}{\overline{z}})}$.
   \end{lem}
   \begin{lem}{\label{lem5}}
   Let $P(z)$ and $F(z)$ be two polynomials such that $degP(z)\leq degF(z) =  n$. If  $F(z)$ has all zeros in $\overline{\mathbb{D}}$ and $|P(z)| \leq |F(z)|$, for $z\in B(\mathbb{D})$. Then for every real and complex number $\alpha, \beta $ with $|\alpha|\leq 1, |\beta| \leq 1$, $R\geq 1$ and $z \in \mathbb{C\backslash D}$
   \begin{align}{\label{3.4}}
      \nonumber \biggl|\Tilde{\mathbb{S}}_a&[P](Rz)-\alpha \Tilde{\mathbb{S}}_a[P](z)+\beta \biggl\{\biggl(\frac{R+1}{2}\biggr)^n -|\alpha|\biggr\}\Tilde{\mathbb{S}}_a[P](z)\biggl|\\ &\leq \biggl|\Tilde{\mathbb{S}}_a[F](Rz)-\alpha \Tilde{\mathbb{S}}_a[F](z)+\beta \biggl\{\biggl(\frac{R+1}{2}\biggr)^n -|\alpha|\biggr\}\Tilde{\mathbb{S}}_a[F](z)\biggr|.
   \end{align}
   \begin{proof}
       In case $R=1$, we have nothing to prove. Therefore we assume that $R > 1$. By hypothesis, $F(z)$ is a polynomial of degree $n$, having all its zeros in $\overline{\mathbb{D}}$ and $|P(z)|\leq |F(z)|$ for $z\in B(\mathbb{D})$. It follows by Rouche's theorem that for every real and complex number $\lambda$ with $|\lambda| > 1$, the polynomial $H(z)= P(z)-\lambda F(z)$ does not vanish in $ \mathbb{ C\backslash \overline{D}}$.\\
       Applying Lemma \ref{lem1} with $k=1$, for every $R > 1$ 
       \begin{align}{\label{3.5}}
       |H(Rz)| \geq  \biggl(\frac{R+1}{2}\biggr)^n |H(z)| ,\ for\ z \in B(\mathbb{D}).
       \end{align}
       Hence for every real and complex number $\alpha$ with $|\alpha|\leq 1$, we have 
       \begin{align}{\label{3.6}}
            \nonumber|H(Rz)-\alpha H(z)|&\geq |H(Rz)|-|\alpha||H(z)|\\
          & > \biggl\{\biggl(\frac{R+1}{2}\biggr)^n -|\alpha|\biggr\}|H(z)|,\ for \ z \in B(\mathbb{D}).
       \end{align}
       Since $H(Re^{i\theta}) \neq 0$ and $\left(\frac{R+1}{2}\right)^n >1$, hence from (\ref{3.5}), we have 
       \[|H(Re^{i\theta})|>|H(e^{i\theta})|,\ for \ R>1 \ and \ 0\leq \theta \leq 2\pi,\]
       equivalently
       \[|H(Rz)|>|H(z)|, \ for \ z \in B(\mathbb{D}) \ and \ R>1.\]
       Since all the zeros of $H(Rz)$ lie in $\mathbb{D}$, it follows by Rouche's theorem for $|\alpha| \leq 1$ that the polynomial $H(Rz)-\alpha H(z)$ does not vanish in $\mathbb{ C\backslash D}$. Again applying  the Rouche's theorem for $|\beta| \leq 1$, it follows from inequality (\ref{3.6}) that the polynomial
       \[G(z)=H(Rz)-\alpha H(z) + \beta \biggl\{\bigg(\frac{R+1}{2}\biggr)^n-|\alpha|\biggr\}H(z),\]
       has all its zeros in $\mathbb{D}$. Using Lemma \ref{lem2}, all the zeros of $\Tilde{\mathbb{S}}_a[G](z)$ lie in $\mathbb{D}$. Replacing $H(z)$ by $P(z)-\lambda F(z)$ and since $\Tilde{\mathbb{S}}_a$ is linear, it follows that the polynomial
       \begin{align*}
       \Tilde{\mathbb{S}}_a&[P](Rz)-\alpha \Tilde{\mathbb{S}}_a[P](z)+\beta\biggl\{\biggl(\frac{R+1}{2}\biggr)^n -|\alpha|\biggr\}\Tilde{\mathbb{S}}_a[P](z)-\\ &\lambda \biggl[\Tilde{\mathbb{S}}_a[F](Rz)-\alpha \Tilde{\mathbb{S}}_a[F](z)
       +\beta\biggl\{\biggl(\frac{R+1}{2}\biggr)^n -|\alpha|\biggr\} \Tilde{\mathbb{S}}_a[F](z)\biggr]
        \end{align*}
        having no zeros in $\mathbb{ C\backslash D}$.\\
        This implies
        \begin{align*}
            \biggl|& \Tilde{\mathbb{S}}_a[P](Rz)-\alpha \Tilde{\mathbb{S}}_a[P](z)+\beta\biggl\{\biggl(\frac{R+1}{2}\biggr)^n -|\alpha|\biggr\}\Tilde{\mathbb{S}}_a[P](z)\biggr|\\ & \leq \biggl|\Tilde{\mathbb{S}}_a[F](Rz)-\alpha \Tilde{\mathbb{S}}_a[F](z)
       +\beta\biggl\{\biggl(\frac{R+1}{2}\biggr)^n -|\alpha|\biggr\} \Tilde{\mathbb{S}}_a[F](z)\biggr|,\ for\ z \in\mathbb{ C\backslash D}.
        \end{align*}
        If this is not true, then there exists a point $z =z_0 \in \mathbb{ C\backslash D}$, such that
        \begin{align*}
            \biggl|& \Tilde{\mathbb{S}}_a[P](Rz_0)-\alpha \Tilde{\mathbb{S}}_a[P](z_0)+\beta\biggl\{\biggl(\frac{R+1}{2}\biggr)^n -|\alpha|\biggr\}\Tilde{\mathbb{S}}_a[P](z_0)\biggr|\\ & > \biggl|\Tilde{\mathbb{S}}_a[F](Rz_0)-\alpha \Tilde{\mathbb{S}}_a[F](z_0)
       +\beta\biggl\{\biggl(\frac{R+1}{2}\biggr)^n -|\alpha|\biggr\} \Tilde{\mathbb{S}}_a[F](z_0)\biggr|.
        \end{align*}
        Since all the zeros of $F(z)$ lie in $ \overline{\mathbb{D}}$, hence (As in case of $H(z))$ all the zeros of 
        \[\Tilde{\mathbb{S}}_a[F](Rz)-\alpha \Tilde{\mathbb{S}}_a[F](z)
       +\beta\biggl\{\biggl(\frac{R+1}{2}\biggr)^n -|\alpha|\biggr\} \Tilde{\mathbb{S}}_a[F](z)\]
      lie in $\mathbb{D}$, for every real and complex number $\alpha, \beta$ and $|\alpha|\leq 1, |\beta|\leq 1$ and $ R>1$. Therefore
       \[\Tilde{\mathbb{S}}_a[F](Rz_0)-\alpha \Tilde{\mathbb{S}}_a[F](z_0)
       +\beta\biggl\{\biggl(\frac{R+1}{2}\biggr)^n -|\alpha|\biggr\} \Tilde{\mathbb{S}}_a[F](z_0) \neq 0\]
       with $ z_0 \in  \mathbb{ C\backslash D}$.
       We choose
       \[\lambda =\frac{\Tilde{\mathbb{S}}_a[P](Rz_0)-\alpha \Tilde{\mathbb{S}}_a[P](z_0)
       +\beta\biggl\{\biggl(\frac{R+1}{2}\biggr)^n -|\alpha|\biggr\} \Tilde{\mathbb{S}}_a[P](z_0)}{\Tilde{\mathbb{S}}_a[F](Rz_0)-\alpha \Tilde{\mathbb{S}}_a[F](z_0)
       +\beta\biggl\{\biggl(\frac{R+1}{2}\biggr)^n -|\alpha|\biggr\} \Tilde{\mathbb{S}}_a[F](z_0)},\] 
       so that $|\lambda|>1$. For this value of $\lambda,\ \Tilde{\mathbb{S}}_a[G](z_0) = 0$, for some $z=z_0 \in \mathbb{ C\backslash D}$, which contradicts the fact that all the zeros of $\Tilde{\mathbb{S}}_a[G](z)$ lie in $\mathbb{D}$. This proves the desired result.
   \end{proof}
   \end{lem}
   \begin{lem}{\label{lem6}}
       If $P(z)$ is a polynomial of degree $n$, which does not vanish in $\mathbb{D}$, then for every real and complex number $\alpha, \beta $ with $|\alpha|\leq 1, |\beta| \leq 1$, $R\geq 1$ and $z \in \mathbb{C\backslash D}$
   \begin{align}{\label{3.7}}
      \nonumber \biggl|\Tilde{\mathbb{S}}_a&[P](Rz)-\alpha \Tilde{\mathbb{S}}_a[P](z)+\beta \biggl\{\biggl(\frac{R+1}{2}\biggr)^n -|\alpha|\biggr\}\Tilde{\mathbb{S}}_a[P](z)\biggl|\\ &\leq \biggl|\Tilde{\mathbb{S}}_a[Q](Rz)-\alpha \Tilde{\mathbb{S}}_a[Q](z)+\beta \biggl\{\biggl(\frac{R+1}{2}\biggr)^n -|\alpha|\biggr\}\Tilde{\mathbb{S}}_a[Q](z)\biggr|,
   \end{align}
   where $Q(z)=z^n \overline{P(\frac{1}{\overline{z}})}$.
   \begin{proof}
       Since $P(z)$ is a polynomial of degree $n$, having all its zeros in $\mathbb{C\backslash D}$. Therefore, all the zeros of the polynomial $Q(z)=z^n \overline{P(\frac{1}{\overline{z}})}$ lie in $\overline{\mathbb{D}}$ and $|P(z)|= |Q(z)|$ for $z \in B(\mathbb{D})$. Applying Lemma \ref{lem5} with $F(z)$ replaced by $Q(z)$, we get the desired result.
   \end{proof}
    \end{lem}
   \begin{center}
   \section{PROOF OF THEOREMS }
 \end{center}
       \begin{proof}[\textbf{Proof of Theorem 2.1}]
       Let $F(z)= Mz^n$, where $M=\max_{z \in B(\mathbb{D})}|P(z)|$, in Lemma \ref{lem5}. We obtain the conclusion of Theorem \ref{thm1}.
     \end{proof}
 \begin{proof}[\textbf{Proof of Theorem 2.2}]
 The result is trivial if $R=1$ (Lemma \ref{lem4}).
So we assume that $R>1$. If \[M =\max_{z\in B(\mathbb{D})}|P(z)|,\] then $|P(z)| \leq M$ for $z \in B(\mathbb{D})$. Now for every real and  complex number $\lambda$ with $|\lambda|>1$, then the polynomial $P(z)+\lambda M$ has no zeros in $\mathbb{D}$ and applying Lemma \ref{lem6}, we have     
\begin{align}{\label{4.3}}
   \nonumber& \biggl|\Tilde{\mathbb{S}}_a[P](Rz)-\alpha\Tilde{\mathbb{S}}_a[P](z)+\beta\biggl\{\bigg(\frac{R+1}{2}\biggr)^n-|\alpha|\biggr\}\Tilde{\mathbb{S}}_a[P](z)-\\ & \nonumber\quad\lambda\biggr[1-\alpha +\beta\bigg\{\biggl(\frac{R+1}{2}\biggr)^n-|\alpha|\bigg\}\bigg]naM\biggr| \\ &\nonumber \leq\biggl|\Tilde{\mathbb{S}}_a[Q](Rz)-\alpha\Tilde{\mathbb{S}}_a[Q](z)+\beta\biggl\{\bigg(\frac{R+1}{2}\biggr)^n-|\alpha|\biggr\}\Tilde{\mathbb{S}}_a[Q](z)+\\ & \quad \overline{\lambda}\bigg[R^n-\alpha+\beta\biggl\{\bigg(\frac{R+1}{2}\bigg)^n-|\alpha|\biggr\}\bigg]M\Tilde{\mathbb{S}}_a[z^n]\biggr|, \ for\ z\in \mathbb{C \backslash D},
\end{align}
 where $|\alpha|\leq 1,|\beta|\leq1$ and $Q(z)=z^n\overline{P(\frac{1}{\overline{z}})}$.
 Choosing the argument of the $\lambda$ on the right-hand side of the above inequality such that
 \begin{align*}
 &\biggl|\Tilde{\mathbb{S}}_a[Q](Rz)-\alpha\Tilde{\mathbb{S}}_a[Q](z)+\beta\biggl\{\bigg(\frac{R+1}{2}\biggr)^n-|\alpha|\biggr\}\Tilde{\mathbb{S}}_a[Q](z)+\\ & \quad \overline{\lambda}\bigg[R^n-\alpha+\beta\biggl\{\bigg(\frac{R+1}{2}\bigg)^n-|\alpha|\biggr\} \biggr]M\Tilde{\mathbb{S}}_a[z^n]\biggr|\\& = |\lambda|\biggl|R^n-\alpha+\beta\biggl\{\bigg(\frac{R+1}{2}\bigg)^n-|\alpha|\biggr\} \biggr|M|\Tilde{\mathbb{S}}_a[z^n]|-\\& \quad\biggl|\Tilde{\mathbb{S}}_a[Q](Rz)-\alpha\Tilde{\mathbb{S}}_a[Q](z)+\beta\biggl\{\bigg(\frac{R+1}{2}\biggr)^n-|\alpha|\biggr\}\Tilde{\mathbb{S}}_a[Q](z)\biggr|.
\end{align*}
From inequality (\ref{4.3}), we get
\begin{align}{\label{4.4}}
    \nonumber& \biggl|\Tilde{\mathbb{S}}_a[P](Rz)-\alpha\Tilde{\mathbb{S}}_a[P](z)+\beta\biggl\{\biggl(\frac{R+1}{2}\biggr)^n-|\alpha|\biggr\}\Tilde{\mathbb{S}}_a[P](z)\biggr|-\\& \nonumber \quad|\lambda|\biggl|1-\alpha +\beta\bigg\{\biggl(\frac{R+1}{2}\biggr)^n-|\alpha|\bigg\}\biggr|n|a|M \\& \nonumber \leq  |\lambda|\biggl|R^n-\alpha+\beta\biggl\{\bigg(\frac{R+1}{2}\bigg)^n-|\alpha|\biggr\}\biggr|M|\Tilde{\mathbb{S}}_a[z^n]|- \\&\qquad\biggl|\Tilde{\mathbb{S}}_a[Q](Rz)-\alpha\Tilde{\mathbb{S}}_a[Q](z)+\beta\biggl\{\bigg(\frac{R+1}{2}\biggr)^n-|\alpha|\biggr\}\Tilde{\mathbb{S}}_a[Q](z)\biggr|,
\end{align}
for $z \in \mathbb{C\backslash D}, |\alpha|\leq1,|\beta| \leq 1$ and $R>1$, taking $|\lambda|\to 1$ inequality (\ref{4.4}), we get the desired result.
 \end{proof}
 \begin{proof}[\textbf{Proof of Theorem 2.3}]
     Since $P(z)$ has no zeros in $\mathbb{D}$, therefore by Lemma \ref{lem6}, we  have
     \begin{align*}
         \biggl|&\Tilde{\mathbb{S}}_a[P](Rz)-\alpha \Tilde{\mathbb{S}}_a[P](z)+\beta \biggl\{\biggl(\frac{R+1}{2}\biggr)^n -|\alpha|\biggr\}\Tilde{\mathbb{S}}_a[P](z)\biggr|\\ &\leq \biggr|\Tilde{\mathbb{S}}_a[Q](Rz)-\alpha \Tilde{\mathbb{S}}_a[Q](z)+\beta \biggl\{\biggl(\frac{R+1}{2}\biggr)^n -|\alpha|\biggr\}\Tilde{\mathbb{S}}_a[Q](z)\biggr|,
     \end{align*}
     for $ z \in \mathbb{C\backslash D}$  and $Q(z)=z^n\overline{P(\frac{1}{\overline{z}})}$.\\
     Equivalently
     \begin{align*}
         2\biggl|&\Tilde{\mathbb{S}}_a[P](Rz)-\alpha \Tilde{\mathbb{S}}_a[P](z)+\beta \biggl\{\biggl(\frac{R+1}{2}\biggr)^n -|\alpha|\biggr\}\Tilde{\mathbb{S}}_a[P](z)\biggr|\\ &\leq \biggl|\Tilde{\mathbb{S}}_a[Q](Rz)-\alpha \Tilde{\mathbb{S}}_a[Q](z)+\beta \biggl\{\biggl(\frac{R+1}{2}\biggr)^n -|\alpha|\biggr\}\Tilde{\mathbb{S}}_a[Q](z)\biggr|+\\& \qquad \biggl|\Tilde{\mathbb{S}}_a[P](Rz)-\alpha \Tilde{\mathbb{S}}_a[P](z)+ \beta \biggl\{\biggl(\frac{R+1}{2}\biggr)^n -|\alpha|\biggr\}\Tilde{\mathbb{S}}_a[P](z)\biggr|.
     \end{align*}
     Appling Theorem \ref{thm2}, we have
     \begin{align*}
        2\biggl|&\Tilde{\mathbb{S}}_a[P](Rz)-\alpha \Tilde{\mathbb{S}}_a[P](z)+\beta \biggl\{\biggl(\frac{R+1}{2}\biggr)^n -|\alpha|\biggr\}\Tilde{\mathbb{S}}_a[P](z)\biggr|\\& \leq\biggl[ \biggl|R^n-\alpha+\beta\biggl\{\bigg(\frac{R+1}{2}\bigg)^n-|\alpha|\biggr\}\biggr||\Tilde{\mathbb{S}}_a[z^n]|+\\& \qquad\biggl|1-\alpha +\beta\bigg\{\biggl(\frac{R+1}{2}\biggr)^n-|\alpha|\bigg\}\biggr|n|a|\biggr]\max_{z \in \mathbb{D}}|P(z)|.
     \end{align*}
     This is the complete proof of Theorem \ref{thm3}. 
 \end{proof}
 \bibliographystyle{amsplain}
 \bibliography{Derivative_Smirnov}
\end{document}